\documentclass[reqno]{amsart}
\usepackage[left=1in, right=1in, top=1in]{geometry}

\usepackage{etoolbox}

\makeatletter
\let\old@tocline\@tocline
\let\section@tocline\@tocline
\newcommand{\section@dotsep}{4.5}
\newcommand{\subsection@dotsep}{4.5}
\patchcmd{\@tocline}
  {\hfil}
  {\nobreak
     \leaders\hbox{$\m@th
        \mkern \section@dotsep mu\hbox{.}\mkern \section@dotsep mu$}\hfill
     \nobreak}{}{}
\let\section@tocline\@tocline
\let\@tocline\old@tocline

\patchcmd{\@tocline}
  {\hfil}
  {\nobreak
     \leaders\hbox{$\m@th
        \mkern \subsection@dotsep mu\hbox{.}\mkern \subsection@dotsep mu$}\hfill
     \nobreak}{}{}
\let\subsection@tocline\@tocline
\let\@tocline\old@tocline

\let\old@l@section\l@section
\let\old@l@subsection\l@subsection

\def\@tocwriteb#1#2#3{%
  \begingroup
    \@xp\def\csname #2@tocline\endcsname##1##2##3##4##5##6{%
      \ifnum##1>\c@tocdepth
      \else \sbox\z@{##5\let\indentlabel\@tochangmeasure##6}\fi}%
    \csname l@#2\endcsname{#1{\csname#2name\endcsname}{\@secnumber}{}}%
  \endgroup
  \addcontentsline{toc}{#2}%
    {\protect#1{\csname#2name\endcsname}{\@secnumber}{#3}}}%

\newlength{\@tocsectionindent}
\newlength{\@tocsubsectionindent}
\newlength{\@tocsubsubsectionindent}
\newlength{\@tocsectionnumwidth}
\newlength{\@tocsubsectionnumwidth}
\newlength{\@tocsubsubsectionnumwidth}
\newcommand{\settocsectionnumwidth}[1]{\setlength{\@tocsectionnumwidth}{#1}}
\newcommand{\settocsubsectionnumwidth}[1]{\setlength{\@tocsubsectionnumwidth}{#1}}
\newcommand{\settocsubsubsectionnumwidth}[1]{\setlength{\@tocsubsubsectionnumwidth}{#1}}
\newcommand{\settocsectionindent}[1]{\setlength{\@tocsectionindent}{#1}}
\newcommand{\settocsubsectionindent}[1]{\setlength{\@tocsubsectionindent}{#1}}
\newcommand{\settocsubsubsectionindent}[1]{\setlength{\@tocsubsubsectionindent}{#1}}

\renewcommand{\l@section}{\section@tocline{1}{\@tocsectionvskip}{\@tocsectionindent}{}{\@tocsectionformat}}%
\renewcommand{\l@subsection}{\subsection@tocline{1}{\@tocsubsectionvskip}{\@tocsubsectionindent}{}{\@tocsubsectionformat}}%
\renewcommand{\l@subsubsection}{\subsubsection@tocline{1}{\@tocsubsubsectionvskip}{\@tocsubsubsectionindent}{}{\@tocsubsubsectionformat}}%
\newcommand{\@tocsectionformat}{}
\newcommand{\@tocsubsectionformat}{}
\newcommand{\@tocsubsubsectionformat}{}
\expandafter\def\csname toc@1format\endcsname{\@tocsectionformat}
\expandafter\def\csname toc@2format\endcsname{\@tocsubsectionformat}
\expandafter\def\csname toc@3format\endcsname{\@tocsubsubsectionformat}
\newcommand{\settocsectionformat}[1]{\renewcommand{\@tocsectionformat}{#1}}
\newcommand{\settocsubsectionformat}[1]{\renewcommand{\@tocsubsectionformat}{#1}}
\newcommand{\settocsubsubsectionformat}[1]{\renewcommand{\@tocsubsubsectionformat}{#1}}
\newlength{\@tocsectionvskip}
\newcommand{\settocsectionvskip}[1]{\setlength{\@tocsectionvskip}{#1}}
\newlength{\@tocsubsectionvskip}
\newcommand{\settocsubsectionvskip}[1]{\setlength{\@tocsubsectionvskip}{#1}}
\newlength{\@tocsubsubsectionvskip}
\newcommand{\settocsubsubsectionvskip}[1]{\setlength{\@tocsubsubsectionvskip}{#1}}

\patchcmd{\tocsection}{\indentlabel}{\makebox[\@tocsectionnumwidth][l]}{}{}
\patchcmd{\tocsubsection}{\indentlabel}{\makebox[\@tocsubsectionnumwidth][l]}{}{}
\patchcmd{\tocsubsubsection}{\indentlabel}{\makebox[\@tocsubsubsectionnumwidth][l]}{}{}

\newcommand{\@sectypepnumformat}{}
\renewcommand{\contentsline}[1]{%
  \expandafter\let\expandafter\@sectypepnumformat\csname @toc#1pnumformat\endcsname%
  \csname l@#1\endcsname}
\newcommand{\@tocsectionpnumformat}{}
\newcommand{\@tocsubsectionpnumformat}{}
\newcommand{\@tocsubsubsectionpnumformat}{}
\newcommand{\setsectionpnumformat}[1]{\renewcommand{\@tocsectionpnumformat}{#1}}
\newcommand{\setsubsectionpnumformat}[1]{\renewcommand{\@tocsubsectionpnumformat}{#1}}
\newcommand{\setsubsubsectionpnumformat}[1]{\renewcommand{\@tocsubsubsectionpnumformat}{#1}}
\renewcommand{\@tocpagenum}[1]{%
  \hfill {\mdseries\@sectypepnumformat #1}}
  
\let\oldappendix\appendix
\renewcommand{\appendix}{%
  \leavevmode\oldappendix%
  \addtocontents{toc}{%
    \protect\settowidth{\protect\@tocsectionnumwidth}{\protect\@tocsectionformat\sectionname\space}%
    \protect\addtolength{\protect\@tocsectionnumwidth}{2em}}%
}
\makeatother

\makeatletter
\settocsectionnumwidth{2em}
\settocsubsectionnumwidth{2.5em}
\settocsubsubsectionnumwidth{3em}
\settocsectionindent{1pc}%
\settocsubsectionindent{\dimexpr\@tocsectionindent+\@tocsectionnumwidth}%
\settocsubsubsectionindent{\dimexpr\@tocsubsectionindent+\@tocsubsectionnumwidth}%
\makeatother

\settocsectionvskip{0pt}
\settocsubsectionvskip{0pt}
\settocsubsubsectionvskip{0pt}

\usepackage{bbding}
\usepackage{mathrsfs}
\usepackage{mathrsfs}
\usepackage{amsfonts}
\usepackage{comment}
\usepackage{cases}
\usepackage{latexsym}
\usepackage{amsmath}
\usepackage[all]{xy}
\usepackage{stmaryrd}
\usepackage{amsfonts}
\usepackage{amsmath,amssymb,amscd,bbm,amsthm,mathrsfs,dsfont}
\usepackage[draft, notcite]{showkeys}

\usepackage{tikz}

\usepackage{pgflibraryarrows}
\usepackage{pgflibrarysnakes}

\usepackage[numbers,sort&compress]{natbib}
\usepackage[pagebackref]{hyperref}
\usepackage{hypernat}
 
\usepackage{color, hyperref}
\definecolor{blue}{rgb}{0.9,0.0,0.9}
\hypersetup{colorlinks, breaklinks,
            linkcolor=red,urlcolor=blue,
            anchorcolor=blue,citecolor=blue}

\usepackage{fancyhdr}
\usepackage{amsxtra,ifthen}
\usepackage{verbatim}

\usepackage{tikz,xcolor,hyperref}
\definecolor{lime}{HTML}{A6CE39}
\DeclareRobustCommand{\orcidicon}{
\begin{tikzpicture}
\draw[lime, fill=lime] (0,0)
circle[radius=0.16]
node[white]{{\fontfamily{qag}\selectfont \tiny \.{I}D}}; 
\end{tikzpicture}
\hspace{-2mm}
}
\foreach \x in {A, ..., Z}{%
\expandafter\xdef\csname orcid\x\endcsname{\noexpand\href{https://orcid.org/\csname orcidauthor\x\endcsname}{\noexpand\orcidicon}}
}

\usepackage{enumitem}

\usepackage{fancyhdr}
\usepackage{amsxtra,ifthen}
\usepackage{verbatim}

\numberwithin{equation}{section}

\theoremstyle{plain}
\newtheorem{theorem}{Theorem}[section]
\newtheorem{lemma}[theorem]{Lemma}
\newtheorem{proposition}[theorem]{Proposition}

\newtheorem{corollary}[theorem]{Corollary}

\theoremstyle{definition}

\newtheorem{remark}[theorem]{Remark}

\allowdisplaybreaks[4]

\newcommand{\Z}{\mathrm{Z}}
\newcommand{\U}{\mathcal{U}}

\begin{document}

\title[Lie $n$-centralizers of von Neumann algebras]
{Lie $n$-centralizers of von Neumann algebras}

\author{Mohammad Ashraf \hspace{-1.5mm}}
\address{Ashraf: Department of Mathematics, Aligarh Muslim University, Aligarh 202002, India}
\email{mashraf80@hotmail.com}

\author{Mohammad Afajal Ansari* \hspace{-1.5mm}}\thanks{*Corresponding Author}
\address{Ansari (Corresponding Author): Department of Mathematics, Aligarh Muslim University, Aligarh 202002, India}
\email{afzalgh1786@gmail.com}

\author{Md Shamim Akhter\hspace{-1.5mm}}
\address{Akhter: Department of Mathematics, Aligarh Muslim University, Aligarh 202002, India}
\email{akhter2805@gmail.com}

\author{Feng Wei\hspace{-1.5mm}}
\address{Wei : School of Mathematics and Statistics, Beijing
Institute of Technology, Beijing, 100081, China}
\email{daoshuo@hotmail.com} \email{daoshuowei@gmail.com}


\date{\today}

\maketitle

 \begin{abstract}
Let $\U$ be a von Neumann algebra with a projection $P\in \U$. For any $A_1,A_2,\ldots,A_n\in\U,$ define $p_1(A_1)=A_1,$ 
$p_n (A_1,A_2,\ldots,A_n)=[p_{n-1} (A_1,A_2,\ldots,A_{n-1}),A_n]$ for all integers $n\geq 2,$ where $[A,B]=AB-BA$ $(A,B\in\U)$ denotes the usual Lie product.
Assume that $\phi:\U\to\U$ is an additive mapping satisfying
\[\phi(p_n(A_1, A_2, \ldots, A_n)) = p_n(\phi(A_1), A_2, \ldots, A_n) = p_n(A_1, \phi(A_2), \ldots, A_n) \]
for all $A_1, A_2, \ldots, A_n \in \U$ with $A_1A_2=P$
In this article, it is shown that the map $\phi$ is of the form $\phi(A)=WA+\xi(A)$ for all $A\in \U$, where $W\in \mathrm{Z}(\U)$, and $\xi:\U \to \Z(\U)$ ($\Z(\U)$ is the center of $\U$) is an additive map such that $\xi(p_n(A_1, A_2, \ldots, A_n) )=0$ for any $A_1, A_2, \ldots, A_n \in \U$ with $A_1A_2=P$. 
As an application, we characterize generalized Lie $n$-derivations on arbitrary von Neumann algebras.
\end{abstract}

\section{\bf Introduction and Preliminaries}
Let $\mathcal{U}$ be an associative algebra with center $\mathcal{Z(U)}.$ 
A linear mapping $\delta:\mathcal{U}\rightarrow \mathcal{U}$  is called a \textit{centralizer} if
$\delta(xy)=\delta(x)y=x\delta(y)$ holds for all $x, y\in \mathcal{U}.$
For any $x_1,x_2,\ldots, x_n\in\mathcal{U},$ define $p_1(x_1)= x_1,$ $p_2(x_1,x_2)=[x_1,x_2]$ and
$p_n(x_1, x_2,\ldots,x_n)=[p_{n-1}(x_1,x_2,\ldots,x_{n-1}),x_n]$ for all integers $n\geq{2},$ where $[x_1,x_2]=x_1x_2-x_2x_1$ denotes the Lie product of $x_1$ and $x_2$ in $\mathcal{U}$. 
A linear mapping $\delta:\mathcal{U}\rightarrow \mathcal{U}$  is said to be a \textit{Lie $n$-centralizer} ($n\geq{2}$) if
$\delta(p_{n}(x_{1},x_{2},\ldots,x_{n}))=p_{n}(\delta(x_{1}),x_{2},\ldots,x_{n})$ holds for all $x_{1},x_{2},\ldots,x_{n}\in\mathcal{U}.$
In particular, a Lie $2$-centralizer is called a Lie centralizer and a Lie $3$-centralizer is said to be a Lie triple centralizer.  
It is easy to see that every centralizer is a Lie $n$-centralizer but the converse is not necessarily true. 
Describing Lie $n$-centralizers in terms of centralizers is an interesting topics of research. 
Fo\v{s}ner and Jing in \cite{fo} have studied non-additive Lie centralizers on triangular rings, and in \cite{liu2} non-linear Lie centralizers on generalized matrix algebra have been checked. Jabeen in \cite{jab} has described the structure of linear Lie centralizers on a generalized matrix algebra under some conditions. In \cite{fad4} it has been shown that under some conditions on a unital generalized matrix algebra $\mathcal{U}$, if $\phi:\mathcal{U}\rightarrow \mathcal{U}$ is a linear Lie triple centralizer, then $\phi(a)=\lambda a+\xi (a) $ in which $\lambda\in \mathrm{Z}(\mathcal{U})$ and $\xi$ is a linear map from $\mathcal{U}$ into $\mathrm{Z}(\mathcal{U})$ vanishing at every second commutator $[[a,b],c]$ for all$a,b,c\in \U$, where $\mathrm{Z}(\mathcal{U})$ is the center of $\U$. In \cite{yu}, Lie n-centralizers of generalized matrix algebras have been examined. Fadaee et al. in \cite{beh} have studied the characterization of Lie centralizers on non-unital triangular algebras through zero products. In \cite{gh3}, linear Lie centralizers at the zero products on some operator algebras are studied, and in \cite{fad3}, linear Lie centralizers through zero products on a $2$-torsion free unital generalized matrix algebra under some mild conditions, are described.  In \cite{fad2}, the problem of characterizing linear maps behaving like Lie centralizers at idempotent products on triangular algebras is considered, and in \cite{goo}, additive Lie centralizers through idempotent-products on a 2-torsion free unital prime ring are determined.

Other important classes of mappings on algebras are derivations, Lie $n$-derivations and their generalizations, which have been extensively studied. 
A linear mapping $\delta:\mathcal{U}\rightarrow \mathcal{U}$  is called a \textit{derivation} if
$\delta(xy)=\delta(x)y+x\delta(y)$ holds for all $x, y\in \mathcal{U}.$
A linear mapping $\delta:\mathcal{U}\rightarrow\mathcal{U}$ is said to be a \textit{Lie $n$-derivation} ($n\geq{2}$) if
$\delta(p_{n}(x_{1},x_{2},\ldots,x_{n})) = p_{n}(\delta(x_{1}),x_{2},\ldots,x_{n})+p_{n}(x_1, \delta(x_{2}),\ldots,x_{n})
  +\cdots +p_{n}(x_1,x_{2},\ldots,\delta(x_{n}))$ holds for all $x_{1},x_{2},\ldots,x_{n}\in\mathcal{U}.$
Lie $n$-derivations have been further generalized as follows:
Let $\xi:\mathcal{U}\rightarrow\mathcal{U}$ be a linear mapping and $\delta$ be a Lie derivation on $\mathcal{U}$. 
Then $\xi$ is called a {\it generalized Lie $n$-derivation} associated with the Lie derivation $\delta$ if 
$\xi(p_{n}(x_{1},x_{2},\ldots,x_{n})) = p_{n}(\xi(x_{1}),x_{2},\ldots,x_{n})+p_{n}(x_1, \delta(x_{2}),\ldots,x_{n})
  +\cdots +p_{n}(x_1,x_{2},\ldots,\delta(x_{n}))$ holds for all $x_{1},x_{2},\ldots,x_{n}\in\mathcal{U}.$
Note that Lie $n$-centralizers and Lie $n$-derivations are special examples of generalized Lie $n$-derivations. 
Determining the structure of generalized Lie $n$-derivations is of special interest. We refer the reader to
\cite{AA21,ash,be12,B18,B19,fwx13,FengQi,WW13,W14} and references therein about characterization of Lie $n$-derivations and generalized
Lie $n$-derivations. Benkovi\v{c} \cite{B18} recently described the form of generalized Lie $n$-derivations through the structure of Lie $n$-centralizers of a unital algebra $\mathcal{U}$ with a nontrivial idempotent  and proved that under certain assumptions every generalized Lie $n$-derivation $G:\mathcal{U}\to\mathcal{U}$ is of the form $G(T)=ZT+D(T)$ for all $T\in\mathcal{U},$ where $Z\in\mathcal{Z}(\mathcal{U})$ and $D:\mathcal{U}\to\mathcal{U}$ is a Lie $n$-derivation.
Motivated by these developments, in the present article, we study Lie $n$-centralizers on arbitrary von Neumann algebras. It should be noted that most of the previous results about von Neumann algebras are on factor von Neumann algebras or von Neumann algebras without central summands of type $ I_1 $, but our results are on a wider class of von Neumann algebras, and some of our results are generalizations of some previous results. It is worth to mention that by using the obtained results, it is possible to characterize generalized Lie $n$-derivations, generalized Lie derivations, Jordan centralizers and Jordan generalized derivations on von Neumann algebras.

A von Neumann algebra $\U$ is a weakly closed, self-adjoint algebra of operators on a complex Hilbert space $ \mathcal {H} $ containing the identity operator $I$. If $P\in \U$ is idempotent (i. e. $P^2 =P$) and self-adjoint ($P^{*}=P$), then $P$ is called a \textit{projection}. A projection $P\in \U$ is said to be a \textit{central abelian projection} if $ P \in  \mathrm{Z} (\U) $ and $ P\U P $ is abelian. The \textit{central carrier} of $T\in \U$ is the smallest central projection $P$ satisfying $ PT = T$, and denoted by $ \overline{T}$. It is well known that $ \overline{T}$ is the projection whose range is the closed linear span of $ \{ AT(h) : A\in \U, h \in   \mathcal {H} \}$. For each self-adjoint operator $ S \in \U$, the \textit{core} of $S$, denoted by $ \underline{S} $, is $ \sup \{W \in  \mathrm{Z} (\U) :  W= W^*, W \leq S \} $. The projection $P$ is a \textit{core-free projection}, if  $ P \in \U$ is a projection and $ \underline{P} = 0$. A routine verifications shows that $\underline{P} = 0$ if and only if $ \overline{I-P} = I$. Note that $\U$ is a von Neumann algebra with no central summands of type $I_1$ if and only if it has a projection $ P$ such that  $\underline{P} = 0$ and $ \overline{P} = I $. If $\U$ is an arbitrary von Neumann algebra, the unit element $I$ of $\U$ is the sum of two orthogonal central projections $E_1$ and $E_2$ such that $\U= \U E_1\oplus \U E_2 $, $ \U E_1$ is of type $I_1$ and $\U E_2 $ is a von Neumann algebra with no central summands of type $I_1$. So $ \U E_2$ contains a core-free projection with central carrier $E_2$. We refer the reader to \cite{ka} for the theory of von Neumann algebras. 
\begin{remark}\label{von}
Let $\U$ be a von Neumann algebra with no central summands of type $I_1$, and $P\in \U$ be a projection such that $\underline{P} = 0$ and $ \overline{P} = I $. We have $\underline{I-P} = 0$ and $ \overline{I-P} = I $. 
\begin{itemize}
\item[(i)] By \cite[Corollary 5.5.7]{ka} we have 
\[\mathrm{Z}(P\U P)=P\mathrm{Z} (\U) \quad \text{and} \quad  \mathrm{Z}((I-P)\U (I-P))=(I-P)\mathrm{Z}(\U).\]
\item[(ii)]  It follows from the definition of the central carrier that both $ span\{AP (h) : A\in \U , h \in \mathcal{H} \}$ and $ span\{A (I-P) (h) : A\in \U, h \in \mathcal{H} \}$ are dense in $\mathcal{H}$. So  $ A \in \U $,  $ A \U P  = \{ 0\}$  implies $ A = 0$ and $ A \U (I-P )   =\{0 \} $  implies $ A  = 0$.
\end{itemize}
\end{remark}

\begin{remark}\label{rem2.1} Let $\U$ be a von Neumann algebra with a projection $P_1\in\U,$ and write $P_2=I-P_1$. Then, for any $A\in\U$, we have
\begin{itemize}
  \item [$(i)$]   $p_n(A, P_1, P_1, \ldots, P_1)=(-1)^{n-1}P_1AP_2+P_2AP_1$;
  \item [$(ii)$]  $p_n(A, P_2, P_2, \ldots, P_2)=P_1AP_2+(-1)^{n-1}P_2AP_1.$
\end{itemize}
\end{remark}

\section{\bf Lie $n$-centralizers of von Neumann algebras}
The following theorem is the main result of this article.
\begin{theorem}\label{main}
Let $\U$ be a von Neumann algebra with unit element $I$, and $E_1 + E_2  =I$, where $E_1$ and $E_2$ are two orthogonal central projections such that $ \U E_1$ is of type $I_1$ and $ \U E_2$ is a von Neumann algebra with no central summands of type $I_1$. Suppose that $P\in  \U E_2$ is a core-free projection with central carrier $E_2$. Let $\phi: \U \rightarrow \U$ be an additive map. Then $\phi$ satisfies
\[\phi(p_n(A_1, A_2, \ldots, A_n)) = p_n(\phi(A_1), A_2, \ldots, A_n) = p_n(A_1, \phi(A_2), \ldots, A_n) \quad (\textbf{P})\]
for all $A_1, A_2, \ldots, A_n \in \U$ with $A_1A_2=P$ if and only if $\phi(A)=WA+\xi(A)$ $(A\in \U)$, where $W\in \mathrm{Z}(\U)$, $\xi:\U \to \mathrm{Z}(\U)$ is an additive map in which $\xi(p_n(A_1, A_2, \ldots, A_n))=0$ for any $A_1, A_2, \ldots, A_n \in \U$ with $A_1A_2=P.$ 
\end{theorem}
The following result is obtained from Theorem \ref{main} which is a generalization of the obtained result in \cite[Remark 4.4]{fad4}, \cite[Corollary 5.2-(iv)]{fad3} and \cite[Corollary 4.3-(ii)]{fos2} for factor von Neumann algebras.
\begin{corollary}\label{cor1 main}
Let $ \mathcal{U} $ be an arbitrary von Neumann algebra, and $ \phi : \mathcal{U} \to \mathcal{U} $ be an additive map. Then $\phi$ is a Lie $n$-centralizer if and only if there exist an element $ W \in \mathrm{Z}( \mathcal{U} ) $ and an additive map $ \xi : \mathcal{U} \to \mathrm{Z}(\mathcal{U} ) $ such that $ \phi (A) = W A + \xi (A) $ for any $ A \in \mathcal{U} $ and 
$\xi(p_n(A_1, A_2, \ldots, A_n))=0$ for any $A_1, A_2, \ldots, A_n \in \U$.
\end{corollary}

First, we show the main result for von Neumann algebras with no central summands of type $I_1$ in the following proposition.
\begin{proposition}\label{no summand I1}
Let $\U$ be a von Neumann algebra with no central summands of type $I_1$, and $P$ be a core-free projection with central carrier $I$. Let $\phi: \U \rightarrow \U$ be an additive map. Then $\phi$ satisfies $(\textbf{P})$ if and only if $\phi(A)=WA+\xi(A)$ $(A\in \U)$, where $W\in \mathrm{Z}(\U)$, $\xi:\U \to \mathrm{Z}(\U)$ is an additive map in which 
$\xi(p_n(A_1, A_2, \ldots, A_n))=0$ for any $A_1, A_2, \ldots, A_n \in \U$ with $A_1A_2=P$. 
\end{proposition}
\begin{proof}
Assume that $\phi$ satisfies $(\textbf{P})$. Set $P_1 := P$ and $P_2 := I- P_1 $.  By Remark \ref{von} $P_2$ is also core free and $ \overline{P_2} = I $. Set $ \mathcal{U}_{ij} = P_i \mathcal{U} P_j$ ($i, j = 1, 2$), then $ \mathcal{U} = \U_{11}+ \U_{12}+ \U_{21}+ \U_{22}$. So any element $A$ of $ \U $ is of the form $ A = A_{11}+A_{12}+A_{21}+A_{22} $ for some $A_{ij} \in \U_{ij} $ ($ i ,  j = 1 ,  2 $). The continuation of the proof in this case is done through the following lemmas.
\begin{lemma}\label{l1}
$ \phi (I) , \phi (P_1)  \in \mathcal{U}_{11} + \mathcal{U}_{22} $.
\end{lemma}
\begin{proof}
Since $I P_1 = P_1$ and $ P_1 P_1 = P_1 $, we get 
\begin{align*}
  0 &=\phi (p_n(I, P_1, \ldots, P_1)) \\
    &=p_n(\phi(I), P_1, \ldots, P_1)\\
    &=(-1)^{n-1}P_1\phi(I)P_2+P_2\phi(I)P_1
\end{align*}
and 
\begin{align*}
  0 &=\phi (p_n(P_1, P_1, \ldots, P_1)) \\
    &=p_n(\phi(P_1), P_1, \ldots, P_1)\\
    &=(-1)^{n-1}P_1\phi(P_1)P_2+P_2\phi(P_1)P_1.
\end{align*}
 Multiplying the above equations once from left by $P_2,$ and once from right by $P_2$, we arrive at $ P_2 \phi (I) P_1 = P_1 \phi (I) P_2 = 0 $ and $ P_2 \phi (P_1) P_1 = P_1 \phi (P_1) P_2 = 0 $, so $ \phi (I) , \phi (P_1) \in \mathcal{U}_{11} + \mathcal{U}_{22} $.
\end{proof}
\begin{lemma}\label{l2}
$ \phi (\mathcal{U}_{ij} ) \subseteq \mathcal{U}_{ij} $, where $ 1 \leq i \neq j \leq 2 $.
\end{lemma}
\begin{proof}
For any $A_{12} \in\mathcal{U}_{12} $, since $( I+ A_{12} ) P_1 = P_1$, by assumption, we see that
\begin{align*}
  \phi (A_{12}) &= \phi ( p_n(I + A_{12}, P_1, \ldots, P_1)) \\
   &= p_n(I + A_{12}, \phi(P_1),P_1, \ldots, P_1)\\
   &=p_n(I, \phi(P_1),P_1, \ldots, P_1)+p_n(A_{12}, \phi(P_1),P_1, \ldots, P_1)\\
   &=p_n(A_{12}, \phi(P_1),P_1, \ldots, P_1)\\
   &=p_{n-1}([A_{12}, \phi(P_1)],P_1, \ldots, P_1)\\
   &=(-1)^{n-2}P_1[A_{12}, \phi(P_1)]P_2+P_2[A_{12}, \phi(P_1)]P_1\\
   &=(-1)^{n-2}P_1[A_{12}, \phi(P_1)]P_2.
\end{align*}
Multiplying the above equation once from left and right by $ P_1 $, once from left and right by $ P_2 $, and once from left by $ P_2 $ and from right by $ P_1 $, we conclude that
 $ P_1 \phi (A_{12})  P_1 = 0 $,  $ P_2 \phi (A_{12})  P_2 = 0  $ and $ P_2 \phi (A_{12})  P_1 =0 $. 
Therefore $ \phi (A_{12} ) =P_1\phi (A_{12} )P_2 \in \mathcal{U}_{12} $.
\par 
 For any $A_{21} \in\mathcal{U}_{21} $, since $ P_1( I+ A_{21} ) = P_1$, we have
\begin{align*}
 \phi (A_{21}) &= \phi ( p_n(P_1,I + A_{21}, P_1, \ldots, P_1)) \\
  &= p_n(\phi(P_1),I + A_{21}, P_1, \ldots, P_1)\\
  &=p_n(\phi(P_1),I, P_1, \ldots, P_1)+p_n(\phi(P_1),A_{21}, P_1, \ldots, P_2)\\
  &=p_n(\phi(P_1),A_{21}, P_1, \ldots, P_1)\\
  &=p_{n-1}([\phi(P_1),A_{21}], P_1, \ldots, P_1)\\
  &=(-1)^{n-2}P_1[\phi(P_1),A_{21}]P_2+P_2[\phi(P_1),A_{21}]P_1\\
  &=P_2[\phi(P_1),A_{21}]P_1.
\end{align*}
Multiplying the above equation once from left and right by $P_1$, once from left and right by $ P_2 $, and once from left by $ P_2 $ and from right by $ P_1 $, and we arrive at $ \phi (A_{21} ) \in \mathcal{U}_{21} $. 
\end{proof}
\begin{lemma}\label{l3}
$ \phi (\mathcal{U}_{ii} ) \subseteq \mathcal{U}_{11} + \mathcal{U}_{22} $, for $ i \in \{ 1 , 2 \} $.
\end{lemma}
\begin{proof}
For any invertible $A_{11} \in \mathcal{U}_{11} $, since $ A_{11} A_{11}^{-1} = P_1 $, we have
\begin{align*}
  0 &=\phi (p_n(A_{11},A_{11}^{-1}, P_1, \ldots, P_1)) \\
    &=p_n(\phi(A_{11}),A_{11}^{-1}, P_1, \ldots, P_1)\\
    &=p_{n-1}([\phi(A_{11}),A_{11}^{-1}], P_1, \ldots, P_1)\\
    &=(-1)^{n-2}P_1[\phi(A_{11}),A_{11}^{-1}]P_2+P_2[\phi(A_{11}),A_{11}^{-1}]P_1\\
    &=(-1)^{n-1}A_{11}^{-1}\phi(A_{11})P_2+P_2\phi(A_{11})A_{11}^{-1}.
\end{align*}
Write $ \phi (A_{11} ) = \sum_{i,j=1}^2 T_{ij} $. It follows from above equation that $ T_{12}  = T_{21}  = 0 $. Consequently, $ \phi (A_{11} ) \in \mathcal{U}_{11} + \mathcal{U}_{22} $, for any invertible $A_{11} \in \mathcal{U}_{11} $. For any $A_{11} \in \mathcal{U}_{11} $,  we may find a sufficiently big number $m$ such that $ m P_1 - A_{11} $ is invertible. Thus, by the above and Lemma \ref{l1}, we have
\[ \phi (A_{11} ) = m \phi (P_1) - \phi ( m P_1 - A_{11} ) \in \mathcal{U}_{11} + \mathcal{U}_{22} . \]
\par 
For any $B_{22} \in \mathcal{U}_{22} $, write $ \phi (B_{22} ) = \sum_{i,j=1}^2 S_{ij} $. 
By the equation  $ ( P_1 + B_{22} )P_1  = P_1 $ and Lemma \ref{l1}, we have 
\begin{align*}
  0 &=\phi (p_n(P_1 + B_{22}, P_1, \ldots, P_1)) \\
    &= p_n(\phi(P_1 + B_{22}), P_1, \ldots, P_1)\\
    &= p_n(\phi(P_1), P_1, \ldots, P_1)+p_n(\phi(B_{22}), P_1, \ldots, P_1)\\
    &=p_n(\phi(B_{22}), P_1, \ldots, P_1)\\
    &=(-1)^{n-1}P_1\phi(B_{22})P_2+P_2\phi(B_{22})P_1.
\end{align*}
 It follows that $ S_{12} = S_{21} = 0 $. Therefore $ \phi (B_{22} ) \in \mathcal{U}_{11} + \mathcal{U}_{22} $.\\
\end{proof}
\begin{lemma}\label{l4}
For $ i = 1, 2 $, there exists an additive map $h_i : \mathcal{U}_{ii} \to \mathrm{Z}(\mathcal{U} )$ such that $ P_j \phi (A_{ii} ) P_j = h_i (A_{ii} ) P_j $ for any $A_{ii} \in \mathcal{U}_{ii}$, where $ 1 \leq i \neq j \leq 2 $.
\end{lemma}
\begin{proof}
 For any invertible element $A_{11} \in \mathcal{U}_{11} $, and each $B_{22} \in \mathcal{U}_{22} $ we have $ A_{11} ( A_{11}^{-1} + B_{22} ) = P_1 $, and  $ ( A_{11} + B_{22} ) A^{-1}_{11} = P_1 $. First, suppose that $n\geq3.$ For each $ C_{ij} \in \mathcal{U}_{ij}  $, where $ (1\leq i \neq j \leq 2 )$, we get
  \begin{align*}
  0 & = \phi ( p_n( A_{11} , A_{11}^{-1} + B_{22}, C_{21}, P_1, \ldots, P_1) ) \\
  &=p_n( \phi(A_{11}) , A_{11}^{-1} + B_{22}, C_{21}, P_1, \ldots, P_1)\\
  &=p_n( \phi(A_{11}) , A_{11}^{-1}, C_{21}, P_1, \ldots, P_1)+p_n( \phi(A_{11}), B_{22}, C_{21}, P_1, \ldots, P_1)\\
  & = \phi  ( p_n(A_{11}, A_{11}^{-1}, C_{21}, P_1, \ldots, P_1) )  + p_n(\phi(A_{11}), B_{22}, C_{21}, P_1, \ldots, P_1) \\
  &=p_n(\phi(A_{11}), B_{22}, C_{21}, P_1, \ldots, P_1)\\
  &=p_{n-2}([[\phi(A_{11}), B_{22}], C_{21}], P_1, \ldots, P_1)\\
  &=(-1)^{n-3}P_1[ [ \phi ( A_{11} ) , B_{22} ] , C_{21} ]P_2+P_2[ [ \phi ( A_{11} ) , B_{22} ] , C_{21} ]P_1\\
  &=P_2[ [ \phi ( A_{11} ) , B_{22} ] , C_{21} ]P_1
  \end{align*}
   and 
   \begin{align*}
  0 & = \phi ( p_n( A_{11}^{-1}  + B_{22}, A_{11}, C_{12}, P_1, \ldots, P_1) ) \\
  &=p_n( \phi(A_{11}^{-1}  + B_{22}), A_{11}, C_{12}, P_1, \ldots, P_1)\\
  &=p_n( \phi(A_{11}^{-1}) , A_{11}, C_{12}, P_1, \ldots, P_1)+p_n( \phi(B_{22}), A_{11}, C_{12}, P_1, \ldots, P_1)\\
  & = \phi  ( p_n(A_{11}^{-1}, A_{11}, C_{12}, P_1, \ldots, P_1) )  + p_n(\phi(B_{22}), A_{11}, C_{12}, P_1, \ldots, P_1) \\
  &=p_n(\phi(B_{22}), A_{11}, C_{12}, P_1, \ldots, P_1)\\
  &=p_{n-2}([[\phi(B_{22}), A_{11}], C_{12}], P_1, \ldots, P_1)\\
  &=(-1)^{n-3}P_1[ [ \phi ( B_{22} ) , A_{11} ] , C_{12} ]P_2+P_2[ [ \phi ( B_{22} ) , A_{11} ] , C_{12} ]P_1\\
  &=(-1)^{n-3}P_1[ [ \phi ( B_{22} ) , A_{11} ] , C_{12} ]P_2.
  \end{align*}
 Considering above equations, and using Lemma \ref{l3}, we arrive at 
 \begin{equation*}
( P_2 \phi (A_{11} ) P_2 B_{22} - B_{22} P_2 \phi (A_{11} ) P_2 ) C P_1 = 0
\end{equation*}
and
\begin{equation*}
( P_1 \phi ( B_{22} ) P_1 A_{11} - A_{11} P_1 \phi (B_{22} ) P_1 ) C P_2  = 0
\end{equation*}
for any $ C \in \mathcal{U} $ and any invertible element $ A_{11} \in \mathcal{U}_{11} $. For any $A_{11} \in \mathcal{U}_{11} $, 
we may find a sufficiently big number $m$ such that $ m P_1 - A_{11} $ is invertible. So,
\begin{equation*}
( P_2 \phi (A_{11} ) P_2 B_{22} - B_{22} P_2 \phi (A_{11} ) P_2 ) C P_1 = 0
\end{equation*}
and
\begin{equation*}
( P_1 \phi ( B_{22} ) P_1 A_{11} - A_{11} P_1 \phi (B_{22} ) P_1 ) C P_2 = 0
\end{equation*}
for all $ C \in \mathcal{U} $, $A_{11} \in \mathcal{U}_{11} $, and $B_{22} \in \mathcal{U}_{22} $. 
From Remark \ref{von}, we conclude that $ P_2 \phi (A_{11} ) P_2 \in \mathrm{Z}(\mathcal{U}_{22}) $ and $ P_1 \phi ( B_{22} ) P_1 \in \mathrm{Z}(\mathcal{U}_{11}) $, and by fact that $ \mathrm{Z}(\mathcal{U}_{22}) = \mathrm{Z}(\mathcal{U} ) P_2 $ and $ \mathrm{Z}(\mathcal{U}_{11}) = \mathrm{Z}(\mathcal{U} ) P_1 $ we have $ P_2 \phi (A_{11} ) P_2 \in \mathrm{Z}(\mathcal{U} ) P_2 $ and $ P_1 \phi (B_{22} ) P_1 \in \mathrm{Z}(\mathcal{U} ) P_1 $. Therefore, for any $A_{11} \in \mathcal{U}_{11} $ and $ B_{22} \in \mathcal{U}_{22} $, there are $ Z_1 , Z_2 \in \mathrm{Z}(\mathcal{U} ) $ such that 
$P_2 \phi (A_{11} ) P_2 = Z_1 P_2$  and $P_1 \phi (B_{22} ) P_1  = Z_2 P_1.$

Now, suppose that $n=2.$ Then
 \begin{align*}
  0 & = \phi ( [p_n( A_{11} , A_{11}^{-1} + B_{22}]) \\
  &=[\phi(A_{11}) , A_{11}^{-1} + B_{22}]\\
  &=[\phi(A_{11}) , A_{11}^{-1}]+[\phi(A_{11}), B_{22}]\\
  & = \phi  ( [p_n(A_{11}, A_{11}^{-1}])  + [\phi(A_{11}), B_{22}] \\
  &=[\phi(A_{11}), B_{22}]=[P_2\phi(A_{11})P_2, B_{22}]
  \end{align*}
   and 
   \begin{align*}
  0 & = \phi ( [p_n( A_{11}^{-1}  + B_{22}]) \\
  &=[\phi(A_{11}^{-1}  + B_{22}), A_{11}]\\
  &=[\phi(A_{11}^{-1}) , A_{11}]+[\phi(B_{22}), A_{11}]\\
  &= \phi  ( [p_n(A_{11}^{-1}, A_{11}])  + [\phi(B_{22}), A_{11}]\\
  &=[\phi(B_{22}), A_{11}]=[P_1\phi(B_{22})P_1, A_{11}].
 \end{align*}
Again, by Remark \ref{von}, we conclude that $P_2 \phi (A_{11} ) P_2 = Z_1 P_2$  and $P_1 \phi (B_{22} ) P_1  = Z_2 P_1$ for some $ Z_1 , Z_2 \in \mathrm{Z}(\mathcal{U} ).$

So we can define the maps $ h_1 : \mathcal{U}_{11} \to \mathrm{Z}( \mathcal{U} ) $ by $ h_1 (A_{11} )= Z_1 $ for any $A_{11} \in \mathcal{U}_{11} $ and $ h_2 : \mathcal{U}_{22} \to \mathrm{Z}( \mathcal{U} ) $ by $ h_2 (B_{22} )= Z_2 $
 for any $B_{22} \in \mathcal{U}_{22} $. Suppose that $ h_1 (A_{11} ) = Z_1 \in \mathrm{Z} (\mathcal{U} ) $ and $ h_1 (A_{11} ) = Z_1' \in \mathrm{Z} (\mathcal{U} ) $. Then we have $ \phi (A_{11} ) - Z_1 \in \mathcal{U}_{11} $ and $ \phi (A_{11} ) - Z_1' \in \mathcal{U}_{11}$. It follows that $ Z_1' - Z_1 = (\phi (A_{11} ) - Z_1 ) - ( \phi (A_{11} ) - Z_1' ) \in \mathcal{U}_{11} \cap \mathrm{Z} (\mathcal{U} ) = \{ 0 \} $. So $ Z_1 = Z_1' $. In a similar way it is proved that $Z_2$ is unique. By the uniqueness of $ Z_1 $ and $ Z_2 $ the maps $ h_1 $ and $ h_2 $ are well-defined. Moreover, from the uniqueness
of $ Z_1 $ and $ Z_2 $ and additivity of $\phi$ it follows that $h_1$ and $h_2$ are additive. Also,
\[ P_2 \phi (A_{11} ) P_2 = h_1 (A_{11} ) P_2 \quad \text{and} \quad P_1 \phi (B_{22} ) P_1 = h_2 (B_{22} ) P_1 . \]
\end{proof} 
Now, for any $ A = A_{11} + A_{12} + A_{21} + A_{22} \in \mathcal{U} $, we define two additive maps $ \xi : \mathcal{U} \to \mathrm{Z} ( \mathcal{U} ) $ and $ \psi : \mathcal{U} \to \mathcal{U} $ by
\[ \xi ( A ) = h_1 ( A_{11} ) + h_2 ( A_{22} ) \quad \text{and} \quad \psi (A) = \phi (A) - \xi(A) . \]
By Lemmas 1-3, it is clear that $ \psi (\mathcal{U}_{ii} ) \subseteq \mathcal{U}_{ii} $ for $i=1,2$, and $ \psi (\mathcal{U}_{ij} ) = \phi ( \mathcal{U}_{ij} )\subseteq \mathcal{U}_{ij} $ for $ 1 \leq i \neq j \leq 2$.
\begin{lemma}\label{l5}
There is an element $W\in \mathrm{Z}( \mathcal{U}) $ such that $\psi(A)=WA$ for all $A\in \U$.
\end{lemma}
\begin{proof}
We divide the proof into the following steps.
\par 
\textit{Step 1.} The following statements hold:
\begin{enumerate}
\item[(i)]
$ \psi (A_{ii} B_{ij} ) = \psi (A_{ii} ) B_{ij} = A_{ii} \psi (B_{ij} ) $ for all $A_{ii} \in \mathcal{U}_{ii} $ and $ B_{ij} \in \mathcal{U}_{ij}$, where $1 \leq i \neq j \leq 2$;
\item[(ii)]
$ \psi (B_{ij} A_{jj} ) = \psi (B_{ij} ) A_{jj} = B_{ij} \psi (A_{jj} ) $ for all $B_{ij} \in \mathcal{U}_{ij} $ and $ A_{jj} \in \mathcal{U}_{jj}$, where $1 \leq i \neq j \leq 2$
\end{enumerate}
For any invertible element $A_{11} \in \mathcal{U}_{11} $ and any $B_{12} \in \mathcal{U}_{12} $, since $ (A_{11}^{-1} + B_{12}) A_{11} = P_1$, by Lemma \ref{l2}, we have
\begin{align*}
(-1)^{n-1}\psi ( A_{11} B_{12} ) & = \phi ( (-1)^{n-1}A_{11} B_{12} ) \\
& = \phi (p_n(A_{11}^{-1} + B_{12} , A_{11}, P_1,\ldots,P_1)) \\
& = p_n(\phi ( A_{11}^{-1} + B_{12} ) , A_{11}, P_1,\ldots,P_1)  \\
&=p_n(\phi ( A_{11}^{-1}) + \phi(B_{12} ) , A_{11}, P_1,\ldots,P_1)\\
&=p_n(\phi ( A_{11}^{-1}), A_{11}, P_1,\ldots,P_1)+p_n(\phi(B_{12} ) , A_{11}, P_1,\ldots,P_1)\\
&=\phi(p_n(A_{11}^{-1}, A_{11}, P_1,\ldots,P_1))+p_n(\phi(B_{12} ) , A_{11}, P_1,\ldots,P_1)\\
&=p_n(\phi(B_{12} ) , A_{11}, P_1,\ldots,P_1)\\
&=p_{n-1}([\phi(B_{12} ) , A_{11}], P_1,\ldots,P_1)\\
&=(-1)^{n-2}P_1[\phi(B_{12} ) , A_{11}]P_2+P_2[\phi(B_{12} ) , A_{11}]P_1\\
&=(-1)^{n-2}P_1[\phi(B_{12} ) , A_{11}]P_2\\
&=(-1)^{n-1}A_{11} \phi (B_{12} )\\
&= (-1)^{n-1}A_{11} \psi (B_{12} )
\end{align*}
and
\begin{align*}
(-1)^{n-1}\psi ( A_{11} B_{12} ) & = \phi ( (-1)^{n-1}A_{11} B_{12} ) \\
& = \phi (p_n(A_{11}^{-1} + B_{12} , A_{11}, P_1,\ldots,P_1)) \\
& = p_n(A_{11}^{-1} + B_{12}, \phi ( A_{11}), P_1,\ldots,P_1)  \\
&=p_n(\phi (A_{11}^{-1}, \phi (A_{11}), P_1,\ldots,P_1)+p_n(B_{12}, \phi (A_{11}), P_1,\ldots,P_1)\\
&=\phi(p_n(A_{11}^{-1}, A_{11}, P_1,\ldots,P_1))+p_n(B_{12}, \phi(A_{11}), P_1,\ldots,P_1)\\
&=p_n(B_{12}, \phi(A_{11}), P_1,\ldots,P_1)\\
&=p_{n-1}([B_{12}, \phi(A_{11})], P_1,\ldots,P_1)\\
&=(-1)^{n-2}P_1[B_{12}, \phi(A_{11})]P_2+P_2[B_{12}, \phi(A_{11})]P_1\\
&=(-1)^{n-2}P_1[B_{12}, \phi(A_{11})]P_2\\
&=(-1)^{n-1}\phi(A_{11})B_{12}\\
&= (-1)^{n-1}\psi(A_{11})B_{12}.
\end{align*}
 For any $ A_{11} \in \mathcal{U}_{11}$, there exists an integer $m$ such that $m P_1 - A_{11}$ is invertible. Note that $mP_1 $ is also invertible. By above results we have $ \psi ( m P_{1} B_{12} ) = m P_1 \psi (B_{12} ) = m \psi (P_1 ) B_{12} $ and $ \psi ( ( m P_{1} - A_{11}) B_{12} ) = (m P_1 - A_{11} ) \psi (B_{12} ) =\psi ( mP_1 - A_{11} ) B_{12} $. Thus, $ \psi ( A_{11} B_{12} ) = A_{11} \psi (B_{12} ) = \psi (A_{11} ) B_{12} $,
for any $ A_{11} \in \mathcal{U}_{11} $ and any $ B_{22} \in \mathcal{U}_{22}$.
\par 
For any invertible element $A_{11} \in \mathcal{U}_{11} $ and any $B_{21} \in \mathcal{U}_{21} $, since $ A_{11} (A_{11}^{-1} + B_{21}) = P_1 $ and 
$p_n(A_{11}, A_{11}^{-1} + B_{21}, P_1,\ldots,P_1)=-B_{21}A_{11}$, and with the similar arguments as above, it can be checked that
\begin{equation*}
\psi ( B_{21} A_{11} ) = B_{21} \psi (A_{11} ) = \psi (B_{21} ) A_{11}
\end{equation*}
for any $A_{11} \in \mathcal{U}_{11} $ and $B_{21} \in \mathcal{U}_{21} $. For any $A_{22} \in \mathcal{U}_{22} $ and $B_{21} \in \mathcal{U}_{21} $, we have $ (P_1 + A_{22} - A_{22} B_{21} ) ( P_1 + B_{21} ) = P_1 $. By properties of $\psi$ and $\xi$, we see that
\begin{align*}
- \psi (B_{21} ) &=-\phi(B_{21} )\\
& =\phi ( p_n (P_1 + A_{22} - A_{22} B_{21} , P_1 + B_{21},P_1,\ldots, P_1)) \\
& = p_n(\phi ( P_1 ) + \phi ( A_{22} ) - \phi( A_{22} B_{21} ), P_1 + B_{21},P_1,\ldots, P_1) \\
& = p_n(\phi ( P_1 ), P_1,P_1,\ldots, P_1)  + p_n(\phi ( P_1 ),B_{21},P_1,\ldots, P_1)  \\
& \quad+ p_n(\phi ( A_{22} ), P_1,P_1,\ldots, P_1)+p_n(\phi ( A_{22} ),B_{21},P_1,\ldots, P_1)\\
&  \quad-p_n(\phi( A_{22} B_{21} ), P_1,P_1,\ldots, P_1) -p_n(\phi( A_{22} B_{21} ),B_{21},P_1,\ldots, P_1)  \\
& = p_n(\phi ( P_1 ),B_{21},P_1,\ldots, P_1)+p_n(\phi ( A_{22} ),B_{21},P_1,\ldots, P_1)  \\
&  \quad-p_n(\phi( A_{22} B_{21} ), P_1,P_1,\ldots, P_1)\\
& = p_n(\psi ( P_1 )+\xi(P_1),B_{21},P_1,\ldots, P_1)+p_n(\psi ( A_{22} )+\xi ( A_{22} ),B_{21},P_1,\ldots, P_1)  \\
&  \quad-p_n(\psi( A_{22} B_{21} ), P_1,P_1,\ldots, P_1)\\
& = p_n(\psi ( P_1 ),B_{21},P_1,\ldots, P_1)+p_n(\psi ( A_{22} ),B_{21},P_1,\ldots, P_1)  \\
&  \quad-p_n(\psi( A_{22} B_{21} ), P_1,P_1,\ldots, P_1)\\
&=-B_{21}\psi(P_1)+\psi(A_{22})B_{21}-\psi(A_{22}B_{21}).
\end{align*}
Hence, we have
\[ \psi ( A_{22} B_{21} ) = \psi (A_{22} ) B_{21} \]
for all $A_{22} \in \mathcal{U}_{22} $ and $B_{21} \in \mathcal{U}_{21} $, because $\psi (B_{21} )=\psi (B_{21}P_1 )=B_{21}\psi(P_1)$. Also,
\[ \psi ( A_{22} B_{21} ) = \psi ( A_{22} B_{21} P_1 ) = A_{22} B_{21} \psi (P_1) = A_{22} \psi ( B_{21} ) , \]
for any $A_{22} \in \mathcal{U}_{22} $ and $B_{21} \in \mathcal{U}_{21} $.
\par 
 For any $A_{22} \in \mathcal{U}_{22} $ and $B_{12} \in \mathcal{U}_{12} $, since
$ ( P_1 + B_{12} ) (P_1 + A_{22} - B_{12}A_{22}  ) = P_1 $ and $p_n(P_1 + B_{12}, P_1 + A_{22} - B_{12}A_{22},P_2,\ldots,P_2)=-B_{12}$, and with the similar arguments as above, it can be checked that
\[ \psi ( B_{12} A_{22} ) = B_{12} \psi ( A_{22} ) = \psi ( B_{12} ) A_{22} . \]
\par 
 \textit{Step 2.} $ \psi ( A_{ii} B_{ii} ) = \psi (A_{ii} ) B_{ii} = A_{ii} \psi ( B_{ii}) $ for all $ A_{ii}, B_{ii} \in \mathcal{U}_{ii} $, where $ i \in \lbrace1,2\rbrace$.
\\
 For any $ A_{ii} , B_{ii} \in \mathcal{U}_{ii} $ and any $ S_{ij} \in \mathcal{U}_{ij} $ ($1 \leq i \neq j \leq 2$), by Step 1, we have
\[ \psi ( A_{ii} B_{ii} S_{ij} ) = \psi ( A_{ii} B_{ii} ) S_{ij} , \]
and on other hand
\[ \psi ( A_{ii} B_{ii} S_{ij} ) = A_{ii} \psi ( B_{ii} S_{ij} ) = A_{ii} \psi (B_{ii} ) S_{ij} , \]
Comparing the above two equations, we see that $ \psi ( A_{ii} B_{ii} ) S_{ij} = A_{ii} \psi (B_{ii} ) S_{ij} $ holds for all $ S_{ij} \in \mathcal{U}_{ij} $. From Remark \ref{von}, it follows that $ \psi ( A_{ii} B_{ii} ) = A_{ii} \psi (B_{ii} ) $ for any $ A_{ii} , B_{ii} \in \mathcal{U}_{ii} $, where $ i = 1,2 $. Also, for any $ A_{ii} , B_{ii} \in \mathcal{U}_{ii} $ and any $ S_{ji} \in \mathcal{U}_{ji} $ ($1 \leq i \neq j \leq 2$), by Step 1, we get
\[ \psi (S_{ji} A_{ii} B_{ii} ) = S_{ji} \psi ( A_{ii} B_{ii} ) , \]
and on other hand
\[ \psi ( S_{ji} A_{ii} B_{ii} ) = \psi ( S_{ji} A_{ii} ) B_{ii} = S_{ji} \psi ( A_{ii} ) B_{ii} , \]
Comparing the above two equations and by Remark \ref{von}, we see that $ \psi ( A_{ii} B_{ii} ) = \psi ( A_{ii} ) B_{ii} $  for any $ A_{ii} , B_{ii} \in \mathcal{U}_{ii} $, where $ i \in \lbrace1,2\rbrace $. 
\par 
\textit{Step 3.} $ \psi ( A_{ij} B_{ji} ) = \psi (A_{ij} ) B_{ji} = A_{ij} \psi ( B_{ji}) $ for all $ A_{ij} \in \mathcal{U}_{ij} $ and $ B_{ji} \in \mathcal{U}_{ji} $, where $ 1 \leq i \neq j \leq 2 $.
\\
 Assume that $ A_{ij} \in \mathcal{U}_{ij} $ and $ B_{ji} \in \mathcal{U}_{ji} $, $ 1 \leq i \neq j \leq 2 $. It follows from Steps 1 and 2 that
\[ \psi ( A_{ij} B_{ji} ) = \psi ( P_i A_{ij} B_{ji} ) = \psi (P_i ) A_{ij} B_{ji} = \psi ( A_{ij} ) B_{ji} , \]
and
\[ \psi ( A_{ij} B_{ji} ) = \psi ( A_{ij} B_{ji} P_i ) = A_{ij} B_{ji} \psi (P_i ) = A_{ij} \psi ( B_{ji} ) . \]
\par 
\textit{Step 4.} The desired result in Lemma \ref{l5} is valid.
\\
From Steps 1-3 and the fact that each $\U_{ij}$ is a invariant subspace for $\psi$, it follows that 
\[ \psi(AB)=A\psi(B)=\psi(A)B\]
for all $A,B\in \U$. Set $W:=\psi(I)$. So
\[ \psi(A)=\psi(AI)=A\psi(I)=AW\quad \text{and} \quad \psi(A)=\psi(IA)=\psi(I)A=WA\]
for all $A,B\in \U$, and $W\in \mathrm{Z}(\U)$.
\end{proof}
\begin{lemma}\label{l6}
 $ \xi ( p_n(A_1, A_2,\ldots,A_n)) = 0$ for all $ A_1, A_2,\ldots,A_n\in \mathcal{U} $ with $A_1A_2 = P$.
\end{lemma}
\begin{proof}
For any $ A_1, A_2,\ldots,A_n \in \mathcal{U} $ with $A_1A_2 = P $, by Lemma \ref{l5} we have
\begin{align*}
\xi ( p_n(A_1, A_2,\ldots,A_n) ) & = \phi ( p_n(A_1, A_2,\ldots,A_n) ) - \psi (p_n(A_1, A_2,\ldots,A_n)) \\
& =p_n(\phi(A_1), A_2,\ldots,A_n) - \psi (p_n(A_1, A_2,\ldots,A_n) ) \\
&=p_n(\psi(A_1)+\xi(A_1), A_2,\ldots,A_n) - \psi (p_n(A_1, A_2,\ldots,A_n) ) \\
&=p_n(\psi(A_1), A_2,\ldots,A_n) - \psi (p_n(A_1, A_2,\ldots,A_n) ) \\
&=Wp_n(A_1, A_2,\ldots,A_n)-Wp_n(A_1, A_2,\ldots,A_n)\\
& = 0 ,
\end{align*}
since $W\in \mathrm{Z}(\U)$.
\end{proof}
Now, by the definition of $\psi$ and Lemma \ref{l5}, we get $ \phi (A) = W A + \xi (A) $ for any $ A \in \mathcal{U} $, where $ W \in \mathrm{Z}( \mathcal{U} ) $. From Lemma \ref{l6}, it follows that $ \xi : \mathcal{U} \to \mathrm{Z}(\mathcal{U} ) $ is an additive mapping in which $\xi(p_n(A_1, A_2,\ldots,A_n))=0$ for any $A_1, A_2,\ldots,A_n \in \U$ with $A_1A_2=P$. So the desired result is valid.
\par 
The converse is clear.
\end{proof}
Now we are ready to present the proof of the main theorem.\\ \\
\textbf{Proof of Theorem \ref{main}.} Let $\phi$ satisfies $(\textbf{P})$. Suppose that $A,B,Y\in \U$ such that $BYE_2=P$. Put $X:=AE_1 + BE_2 $. Since $\U E_1 \subseteq \Z(\U)$, it follows that $[X,YE_2]=[BE_2,YE_2]=[B,Y]E_2$. Moreover, we have $XYE_2=BYE_2=P.$ Hence, for any $C_3,C_4,\ldots,C_n\in\U,$ we get 
$\phi(p_n(X,YE_2,C_3,C_4,\ldots,C_n))=p_n(\phi(X),YE_2,C_3,C_4,\ldots,C_n)=p_n(X,\phi(YE_2),C_3,C_4,\ldots,C_n)$.
 So
 \begin{align*}
   \phi(p_n(BE_2,YE_2,C_3,C_4,\ldots,C_n)) & =\phi(p_n(X,YE_2,C_3,C_4,\ldots,C_n)) \\
    & =p_n(\phi(X),YE_2,C_3,C_4,\ldots,C_n)\\
    &=p_n(\phi(AE_1),YE_2,C_3,C_4,\ldots,C_n)+p_n(\phi(BE_2),YE_2,C_3,C_4,\ldots,C_n).
 \end{align*}
 Multiplying the right side of the above identity by $E_2$, we obtain
\begin{equation}\label{v1}
\phi(p_n(BE_2,YE_2,C_3,C_4,\ldots,C_n))E_2=p_n(\phi(AE_1),YE_2,C_3,C_4,\ldots,C_n)+p_n(\phi(BE_2),YE_2,C_3,C_4,\ldots,C_n).
\end{equation}
By setting $A=0$ in \eqref{v1} we see that 
\begin{equation}\label{v2}
\phi(p_n(BE_2,YE_2,C_3,C_4,\ldots,C_n))E_2=p_n(\phi(BE_2),YE_2,C_3,C_4,\ldots,C_n).
\end{equation} 
for all $B,Y\in \U E_2 $ with $BYE_2=P$. Also, we have
\begin{equation*}
\phi(p_n(BE_2,YE_2,C_3,C_4,\ldots,C_n))E_2=p_n(AE_1,\phi(YE_2),C_3,C_4,\ldots,C_n)+p_n(BE_2,\phi(YE_2),C_3,C_4,\ldots,C_n).
\end{equation*}
So by the fact that $\U E_1 \subseteq \Z(\U)$ we arrive at
\begin{equation}\label{v3}
\phi(p_n(BE_2,YE_2,C_3,C_4,\ldots,C_n))E_2=p_n(BE_2,\phi(YE_2),C_3,C_4,\ldots,C_n).
\end{equation} 
for all $B,Y\in \U E_2 $ with $BYE_2=P$. Equations \eqref{v2} and \eqref{v3} show that the additive mapping $\varphi: \U E_2\rightarrow \U E_2$ defined by $\varphi(AE_2)=\phi(AE_2)E_2$, on $\U E_2$ satisfies $(\textbf{P})$. By our assumption $ \U E_2$ is a von Neumann algebra with no central summands of type $I_1$, and $P\in \U E_2$ is a projection such that $\underline{P} = 0$ and $ \overline{P} = E_2 $. So by Proposition \ref{no summand I1}, there are $W_1\in \Z(\U E_2)\subseteq \Z(\U)$ and an additive mapping $\xi_1 :\U E_2 \rightarrow \Z(\U E_2)\subseteq \Z(\U)$ such that
\begin{equation}\label{v4}
\phi(AE_2)E_2=\varphi(AE_2)=W_1AE_2+\xi_1(AE_2)
\end{equation}
for all $A\in \U$ and $\xi_1(p_n(AE_2, BE_2,C_3E_2,\ldots,C_nE_2))=0 $ for all $A,B,C_3,\ldots,C_n\in \U$ with $ABE_2=P$. Assume that for $Y\in \U$ the element $YE_2$ is invertible in $\U E_2$ (i.e., $YE_2\in Inv(\mathcal{U} E_2$)). Taking $B:=P(YE_2)^{-1}$ and $X=AE_1+BE_2$ for $A\in \U$. So $BYE_2=P$, and from \eqref{v1} and \eqref{v2}, for any $C_3,\ldots,C_n\in \U$ it follows that $p_n(\phi(AE_1)E_2,YE_2,C_3,\ldots,C_n)=0$. Since  each element of $\U E_2$ is a sum of two invertible elements of $\U E_2$, it results that
$p_n(\phi(AE_1)E_2,YE_2,C_3,\ldots,C_n)=0$
for all $A,Y,C_3,\ldots,C_n\in \U$. However, by Kleinecke-Shirokov Theorem and the fact that the spectral radius is submultiplicative on commuting elements, it results that
$[\phi(AE_1)E_2,YE_2]=0$
for all $A,Y\in \U$. So
\[ \phi(AE_1)E_2\in \Z(\U E_2)\subseteq \Z(\U)\]
for all $A\in \U$. Also 
\[ \varphi(A)E_1\in \U E_1\subseteq \Z(\U)\]
for all $A\in \U$. Now by \eqref{v4} we have
\begin{equation*}
\begin{split}
\phi(A)&=\phi(A)E_1 +\phi(AE_1)E_2 +\phi(AE_2)E_2\\&
=\phi(A)E_1 +\phi(AE_1)E_2 +W_1AE_2+\xi_1(AE_2)\\&
=WA+\xi(A),
\end{split}
\end{equation*}
for all $A\in \U$, where $W:=W_1E_2\in \Z(\U)$ and $\xi: \U\rightarrow \U$ is an additive map defined by $\xi(A)=\phi(A)E_1 +\phi(AE_1)E_2 +\xi_1(AE_2)$. By above all three summands lie in $\Z(\U)$, thus $\xi$ maps $\U$ into $\Z(\U)$. Finally according to these results for $A,B,C_3,\ldots,C_n\in \U$ where $AB=P$ we have
\begin{align*}
\xi(p_n(A,B,C_3,\ldots,C_n))&=\phi(p_n(A,B,C_3,\ldots,C_n))-Wp_n(A,B,C_3,\ldots,C_n)\\
&=p_n(\phi(A),B,C_3,\ldots,C_n)-p_n(WA,B,C_3,\ldots,C_n)\\
&=p_n(\xi(A),B,C_3,\ldots,C_n)=0.
\end{align*}
\par 
The converse is clear. $\hspace{24em}\square$
\\ \\
\textbf{Proof of Corollary \ref{cor1 main}.} Let $\phi$ be a Lie $n$-centralizer. If $\U$ is an abelian von Neumann algebra, then $\phi$ maps $\U$ into $\Z(\U)=\U$. Also, from the fact that $\phi$ is a Lie $n$-centralizer, it follows that $\phi(p_n(A_1, A_2, \ldots, A_n))=0$ for any $A_1, A_2, \ldots, A_n \in \U$. So in this case the result is valid. Now let's assume that $\U$ is non-abelian. In this case the unit element $I$ of $\U$ is the sum of two orthogonal central projections $E_1$ and $E_2$ such that $\U= \U E_1\oplus \U E_2 $, $ \U E_1$ is of type $I_1$ and $\U E_2 $ is a von Neumann algebra with no central summands of type $I_1$. So there exist a core-free projection $P\in  \U E_2$ with central carrier $E_2$. Because $\phi$ is a Lie $n$-centralizer, it satisfies the condition $(\textbf{P})$ on $\U$. Hence by Theorem \ref{main}, $\phi(A)=WA+\xi(A)$ for all $A\in \U$, where $W\in \Z(\U)$, $\xi:\U \to \Z(\U)$ is an additive map. It is sufficient to prove that for any $A_1, A_2, \ldots, A_n\in \U$ we have $\xi(p_n(A_1, A_2, \ldots, A_n))=0$. This part can be proved in similar manner with the proof of Theorem \ref{main}.
\par 
The converse is clear. $\hspace{24em}\square$
\\ \\

\section{\bf Applications}
In this section, we present some applications of the results obtained in the previous section. As an application of Theorem \ref{main}, we characterize the 
generalized Lie $n$-derivations of arbitrary von Neumann algebras.
First, we characterize generalized Lie $n$-derivations of von Neumann algebras with no central summands of type $I_1,$
which is a partial generalization of \cite[Theorem 3.11]{AA21}

\begin{theorem}\label{main2}
Let $\U$ be a von Neumann algebra with no central summands of type $I_1$, and $P$ be a core-free projection with central carrier $I$. Let $\mathcal{L},G_{\mathcal{L}}: \U \rightarrow \U$ be additive maps satisfying
\begin{eqnarray}\label{eqn4.1}
\mathcal{L}(p_n(A_1, A_2, \ldots, A_n))=\sum\limits_{i=1}^{n}p_n(A_1, A_2, \ldots, A_{i-1}, \mathcal{L}(A_i), A_{i+1},\ldots, A_n)
\end{eqnarray}
and
\begin{eqnarray}\label{eqn4.2}
\nonumber G_{\mathcal{L}}(p_n(A_1, A_2, \ldots, A_n))&=&p_n(G_{\mathcal{L}}(A_1), A_2, \ldots, A_n)+\sum\limits_{i=2}^{n}p_n(A_1, A_2, \ldots, A_{i-1}, \mathcal{L}(A_i), A_{i+1},\ldots, A_n) \\
\nonumber &=&p_n(\mathcal{L}(A_1), A_2, \ldots, A_n)+p_n(A_1, G_{\mathcal{L}}(A_2), \ldots, A_n)\\
 &&+\sum\limits_{i=3}^{n}p_n(A_1, A_2, \ldots, A_{i-1}, \mathcal{L}(A_i), A_{i+1},\ldots, A_n)
\end{eqnarray}
for all $A_1, A_2, \ldots, A_n \in \U$ with $A_1A_2=P.$ Then $G_{\mathcal{L}}$ is of the form $G_{\mathcal{L}}(A)=WA+\delta(A)+\chi(A)$ $(A\in \U)$, where $W\in \mathrm{Z}(\U)$, $\delta:\U \to \U$ is an additive derivation and $\chi:\U \to \Z(\U)$ is an additive map satisfying $\chi(p_n(A_1, A_2, \ldots, A_n))=0$ for any $A_1, A_2, \ldots, A_n \in \U$ with $A_1A_2=P.$ 
\end{theorem}
\textbf{Proof.} Let $\phi=G_{\mathcal{L}}-\mathcal{L}.$ Then, $\phi$ satisfies $(\textbf{P})$. Hence, it follows from Proposition \ref{no summand I1} that $\phi$ is of the form
$\phi(A)=WA+\xi(A)$ $(A\in \U)$, where $W\in \mathrm{Z}(\U)$, $\xi:\U \to \mathrm{Z}(\U)$ is an additive map satisfying $\xi(p_n(A_1, A_2, \ldots, A_n))=0$ for any $A_1, A_2, \ldots, A_n \in \U$ with $A_1A_2=P.$ Therefore, $G_{\mathcal{L}}(A)=WA+\xi(A)+\mathcal{L}(A)$ for all $A\in \U.$ By \cite[Theorem 3.11]{AA21}, $\mathcal{L}(A)=\delta(A)+\gamma(A),$ where $\delta:\U \rightarrow \U$ is an additive derivation and $\gamma:\U \rightarrow \mathrm{Z}(\U)$ is an additive map satisfying $\gamma(p_n(A_1, A_2, \ldots, A_n))=0$ for any $A_1, A_2, \ldots, A_n \in \U$ with $A_1A_2=P.$
Suppose that $\chi=\xi+\gamma$. Then $\chi:\U \rightarrow \mathrm{Z}(\U)$ is an additive map satisfying $\chi(p_n(A_1, A_2, \ldots, A_n))=0$ for any $A_1, A_2, \ldots, A_n \in \U$ with $A_1A_2=P$ and $G_{\mathcal{L}}(A)=WA+\delta(A)+\chi(A)$ $(A\in \U),$ as desired.
 $\hspace{24em}\square$

The following corollary is an easy consequence of the above theorem.
\begin{corollary}\label{cor1 main2}
Let $\U$ be a von Neumann algebra with no central summands of type $I_1$, and $ G_{\mathcal{L}} : \mathcal{U} \to \mathcal{U} $ be a generalized Lie $n$-derivation (with associated Lie $n$-derivation $\mathcal{L}$). Then $G_{\mathcal{L}}$ is of the form $G_{\mathcal{L}}(A)=WA+\delta(A)+\chi(A)$ $(A\in \U)$, where $W\in \mathrm{Z}(\U)$, $\delta:\U \to \U$ is an additive derivation and $\chi:\U \to \Z(\U)$ is an additive map satisfying $\chi(p_n(A_1, A_2, \ldots, A_n))=0$ for any $A_1, A_2, \ldots, A_n \in \U.$ 
\end{corollary}

Now, let us state and prove Theorem \ref{main2} for an arbitrary von Neumann algebra.
\begin{theorem}\label{main3}
Let $\U$ be a von Neumann algebra with unit element $I$, and $E_1 + E_2  =I$, where $E_1$ and $E_2$ are two orthogonal central projections such that $ \U E_1$ is of type $I_1$ and $ \U E_2$ is a von Neumann algebra with no central summands of type $I_1$. Suppose that $P\in  \U E_2$ is a core-free projection with central carrier $E_2$. Let $\mathcal{L},G_{\mathcal{L}}: \U \rightarrow \U$ be additive maps satifying \eqref{eqn4.1} and \eqref{eqn4.2}.
Then $G_{\mathcal{L}}$ is of the form $G_{\mathcal{L}}(A)=WA+\chi(A)$ $(A\in \U)$, where $W\in \mathrm{Z}(\U)$, $\chi:\U \to \U$ is an additive map satisfying \eqref{eqn4.1} for any $A_1, A_2, \ldots, A_n \in \U$ with $A_1A_2=P.$ 
\end{theorem}
\textbf{Proof.} Let $\phi=G_{\mathcal{L}}-\mathcal{L}.$ Then, $\phi$ satisfies $(\textbf{P})$. Hence, it follows from Theorem \ref{main} that $\phi$ is of the form
$\phi(A)=WA+\xi(A)$ $(A\in \U)$, where $W\in \mathrm{Z}(\U)$, $\xi:\U \to \mathrm{Z}(\U)$ is an additive map satisfying $\xi(p_n(A_1, A_2, \ldots, A_n))=0$ for any $A_1, A_2, \ldots, A_n \in \U$ with $A_1A_2=P.$ Therefore, $G_{\mathcal{L}}(A)=WA+\xi(A)+\mathcal{L}(A)$ for all $A\in \U.$ Set $\chi=\xi+\mathcal{L}.$ Then it is easy to see that $\chi:\U \to \U$ is an additive map satisfying \eqref{eqn4.1} for any $A_1, A_2, \ldots, A_n \in \U$ with $A_1A_2=P.$ Consequently, $G_{\mathcal{L}}(A)=WA+\chi(A)$ $(A\in \U)$, as desired.
 $\hspace{24em}\square$

\begin{corollary}\label{cor1 main3}
Let $ \mathcal{U} $ be an arbitrary von Neumann algebra,and $ G_{\mathcal{L}} : \mathcal{U} \to \mathcal{U} $ be a generalized Lie $n$-derivation (with associated Lie $n$-derivation $\mathcal{L}$). Then there exist an element 
$ W \in \mathrm{Z}( \mathcal{U} ) $ and a Lie $n$-derivation $ \chi : \mathcal{U} \to \mathcal{U}$ such that $ \phi (A) = W A + \chi (A) $ for any $ A \in \mathcal{U}. $
\end{corollary}
\textbf{Proof.} 
Let $ G_{\mathcal{L}} : \mathcal{U} \to \mathcal{U} $ be a generalized Lie $n$-derivation with associated Lie $n$-derivation $\mathcal{L}.$ 
Then $\phi=G_{\mathcal{L}}-\mathcal{L}$ is a Lie $n$-centralizer on $\U$. Hence, it follows from Corollary \ref{cor1 main2} that there exist an element $ W \in \mathrm{Z}( \mathcal{U} ) $ and an additive map $ \xi : \mathcal{U} \to \mathrm{Z}(\mathcal{U} ) $ such that $ \phi (A) = W A + \xi (A) $ for any $ A \in \mathcal{U} $ and 
$\xi(p_n(A_1, A_2, \ldots, A_n))=0$ for any $A_1, A_2, \ldots, A_n \in \U$. Therefore, $G_{\mathcal{L}}(A)=WA+\xi(A)+\mathcal{L}(A)$ for all $A\in \U.$ 
Set $\chi=\xi+\mathcal{L}.$ Then it is easy to see that $\chi:\U \to \U$ is a Lie $n$-derivation. Consequently, $G_{\mathcal{L}}(A)=WA+\chi(A)$ $(A\in \U)$, as desired.
 $\hspace{24em}\square$

\subsection*{Declarations}
\begin{itemize}
\item[•] All authors contributed to the study conception and design and approved the final manuscript. 
\item[•] The authors have no relevant financial or non-financial interests to disclose.
\item[•] Data sharing not applicable to this article as no datasets were generated or analysed during the current study. 
\end{itemize}

\bibliographystyle{amsplain}

\end{document}